\documentclass[12pt]{amsart}
\usepackage{amsfonts,amscd,amssymb,amsmath,amsthm,mathrsfs,xcolor,lscape,amsmath,amssymb,latexsym}
\usepackage[utf8]{inputenc}
\usepackage{hyperref}
\usepackage[english]{babel}
\newcommand{\Sp}{\mathrm{Sp}}

\newcommand{\GL}{\mathrm{GL}}
\newcommand{\SL}{\mathrm{SL}}

\newcommand{\R}{\mathbb{R}}
\newcommand{\Q}{\mathbb{Q}}
\newcommand{\Z}{\mathbb{Z}}

\newcommand{\C}{\mathbb{C}}
\newcommand{\e}{\epsilon}

\newtheorem{theorem}{Theorem}[section]

\theoremstyle{definition}

\newtheorem{rmk}[theorem]{Remark}

\author{Sandip  Singh}
\address{Department of Mathematics, IIT Bombay, Mumbai}
\email{sandip@math.iitb.ac.in}  

\begin{document}\title[Hypergeometric Monodromy Groups in $\Sp(4)$]{Arithmeticity of some Hypergeometric Monodromy Groups in $\Sp(4)$}

\begin{abstract} The article \cite{SV} gives a list of $51$ symplectic hypergeometric monodromy groups corresponding to primitive pairs of degree four polynomials, which are products of cyclotomic polynomials, and for which, the absolute value of the leading coefficient of the difference polynomial is greater than $2$. 

It follows from \cite{SS} and \cite{SV} that $12$ of the $51$ monodromy groups are arithmetic (cf. Table \ref{table:arithmetic}); and the thinness of $13$ of the remaining $39$ monodromy groups follows from  \cite{BT} (cf. Table \ref{table:thin}). 

In this article, we show that $15$ of the remaining $26$ monodromy groups are arithmetic (cf. Table \ref{table:newarithmetic}). 
\end{abstract}
\maketitle

\section{Introduction} To describe the results obtained in this paper, we first recall some basic facts about the hypergeometric monodromy groups (cf. \cite{BH}).

Let $\alpha=(\alpha_1,\ldots,\alpha_n), \beta=(\beta_1,\ldots,\beta_n)\in\Q^n$, and consider the hypergeometric differential equation
\begin{eqnarray}\label{introdifferentialequation}
D(\alpha;\beta)u=0
\end{eqnarray} on $\mathbb{P}^1(\C)$, where
$D(\alpha;\beta)=(\theta+\beta_1-1)\cdots(\theta+\beta_n-1)-z(\theta+\alpha_1)\cdots(\theta+\alpha_n)$, and $\theta=z\frac{d}{dz}$. 

Note that the differential equation (\ref{introdifferentialequation}) has regular singularities at the points $0, 1$, and $\infty$, and it is regular elsewhere. Therefore, we get a monodromy action of the fundamental group $\pi_1$ of $\mathbb{P}^1(\C)\backslash\{0,1,\infty\}$ on the (local) solution space of the differential equation (\ref{introdifferentialequation}), that is, there is a representation $\rho$ of $\pi_1$ inside $\GL_n(\C)$, which is called the monodromy representation, and its image $\rho(\pi_1)$ is called the monodromy group of the hypergeometric differential equation. Note also that the monodromy group of an $n$-th order hypergeometric differential equation is defined up to conjugation in $\GL_n(\C)$, and it is also called the hypergeometric monodromy group.

The monodromy group $\rho(\pi_1)$ is generated by the monodromy matrices $\rho(h_0)$, $\rho(h_1)$, $\rho(h_\infty)$, where $h_0, h_1, h_\infty$ (loops around $0, 1, \infty$, respectively) are the generators of $\pi_1$, with a single relation $h_\infty h_1 h_0=1$. 

The generators (in $\GL_n(\C$)) of the hypergeometric monodromy group are given by the following theorem of Levelt (\cite{Le}; cf. \cite[Theorem 3.5]{BH}):
\begin{theorem}
 If $\alpha_1,\ldots,\alpha_n,\beta_1,\ldots,\beta_n\in\C$ such that $\alpha_j-\beta_k\not\in\Z$, for all $j,k=1,2,\ldots,n$,  then there exists a basis of the (local) solution space of the hypergeometric differential equation, with respect to which, the actions of $h_\infty$ and $h_0^{-1}$ are, respectively, given by the matrices $A$ and $B$ which are, respectively, the companion matrices of the polynomials  \[f(X)=\prod_{j=1}^{n}(X-{\rm{e}^{2\pi i\alpha_j}})\quad\mbox{ and }\quad g(X)=\prod_{j=1}^{n}(X-{\rm{e}^{2\pi i\beta_j}})\] and the action of $h_1$ is given by $A^{-1}B$.
\end{theorem}

In this article, we consider the cases where $n=4$, and the monodromy groups $\Gamma(f,g)$, with respect to the Levelt's basis, are subgroups of $\SL_4(\Z)$. In particular, let $f,g\in\Z[X]$ be a pair of degree four polynomials, which are products of cyclotomic polynomials, do not have any common root in $\C$, $f(0)=g(0)=1$, and form a primitive pair, that is, $f(X)\neq f_1(X^k)$ and $g(X)\neq g_1(X^k)$, for any $k\geq2$ and $f_1,g_1\in\Z[X]$. 

We now form the companion matrices $A, B$ of $f, g$, respectively, and consider the subgroup $\Gamma(f,g)\subset\SL_4(\Z)$ generated by $A$ and $B$. It follows from Beukers and Heckman \cite[Theorem 6.5]{BH} that $\Gamma(f,g)$ preserves a non-degenerate integral symplectic form $\Omega$ on $\Z^4$ and $\Gamma(f,g)\subset\Sp_4(\Omega)(\Z)$ is Zariski dense in the symplectic group $\Sp_4(\Omega)$ of the form $\Omega$. 

The group $\Gamma(f,g)\subset\Sp_4(\Omega)(\Z)$ is called {\it arithmetic} if it has finite index in $\Sp_4(\Omega)(\Z)$, and {\it thin} if it has infinite index in $\Sp_4(\Omega)(\Z)$. 

In this article, we consider a special case of the question in \cite{S} to determine the pairs of polynomials $f,g$, satisfying the above conditions, and for which, the associated monodromy groups $\Gamma(f,g)$ are arithmetic. The question in \cite{S} has also been considered in the articles \cite{BS}, \cite{BT}, \cite{YYCE}, \cite{E}, \cite{F}, \cite{FMS}, \cite{HvS}, \cite{SS}, \cite{SS0}, \cite{SV}, \cite{Ve2}, \cite{Ve3}.

The article \cite{SV} gives a list \cite[Table 1]{SV} of $60$ pairs of polynomials $f,g$, satisfying the above conditions, and for which, the absolute value of the leading coefficient of the difference $f-g$ is $\leq 2$; and it also gives another list \cite[Table 2]{SV} of $51$ pairs of polynomials $f,g$, satisfying the above conditions, and for which, the absolute value of the leading coefficient of the difference $f-g$ is $\geq 3$. 

The arithmeticity of the monodromy groups associated to the pairs $f,g$ in \cite[Table 1]{SV} follows from \cite[Theorem 1.1]{SV}, which says that, if the absolute value of the leading coefficient of the difference $f-g$ is $\leq2$, then the associated monodromy group $\Gamma(f,g)$ is arithmetic. Note that the pairs $f,g$ in \cite[Table 2]{SV} do not satisfy the condition of \cite[Theorem 1.1]{SV}, as for them, the absolute value of the leading coefficient of the difference $f-g$ is $\geq 3$.

There have been some progress to answer the question to determine the pairs $f,g$ in \cite[Table 2]{SV}, which correspond to an arithmetic or thin $\Gamma(f,g)$. Therefore we split \cite[Table 2]{SV} in $4$ subtables, that is, in Tables \ref{table:arithmetic}, \ref{table:thin}, \ref{table:newarithmetic}, and \ref{table:unknown} of this article. Table \ref{table:arithmetic} lists the pairs $f,g$ of \cite[Table 2]{SV}, for which, the {\it arithmeticity} of the associated monodromy groups follows from \cite{SS} and \cite{SV} (cf. Remark \ref{remarkf(-x)g(-x)}), and Table \ref{table:thin} lists the pairs $f,g$ of \cite[Table 2]{SV}, for which, the {\it thinness} of the associated monodromy groups follows from Brav and Thomas \cite{BT} (cf. Remark \ref{remarkf(-x)g(-x)}). We note here that \cite[Table 2]{SV} contains both $f(X), g(X)$ and $f(-X), g(-X)$, and once we prove the arithmeticity or thinness of $\Gamma(f,g)$ associated to the pairs $f(X), g(X)$, the same follows for $\Gamma(\underbar{f}, \underbar{g})$, where $\underbar{f}, \underbar{g}$ are the polynomials $f(-X), g(-X)$ respectively. For a reference, we give a proof of this observation in the following remark:

\begin{rmk}\label{remarkf(-x)g(-x)}
Let $\underbar{f}, \underbar{g}$ be the polynomials $f(-X)$, $g(-X)$, respectively. Then, we show that the monodromy group $\Gamma(f,g)\subset\Sp_4(\Omega)(\Z)$ is arithmetic if and only if $\Gamma(\underbar{f}, \underbar{g})\subset\Sp_4(^tS\Omega S)(\Z)$ is arithmetic, where $S$ is a $4\times4$ diagonal matrix, which has odd diagonal entries $1$ and even diagonal entries $-1$. This can be shown as follows: Let $A, B, \underbar{A}, \underbar{B}$ be the companion matrices of $f,g, \underbar{f}, \underbar{g}$, respectively, and $C, \underbar{C}$ be the matrices $A^{-1}B$,  $\underbar{A}^{-1}\underbar{B}$, respectively. Then, it can be checked easily that $$S^{-1}(-A)S=\underbar{A},\quad S^{-1}CS=\underbar{C},$$ and the group $<-A,C>$, generated by $-A$ and $C$, is also a subgroup of $\Sp_4(\Omega)(\Z)$. 

Since $\Gamma(f,g)=<A,C>$ has finite index in $\Sp_4(\Omega)(\Z)$ if and only if $<-A,C>$ has finite index in $\Sp_4(\Omega)(\Z)$, and it follows from the above computations that $<-A,C>$ has finite index in $\Sp_4(\Omega)(\Z)$ if and only if $<\underbar{A},\underbar{C}>$ has finite index in $\Sp_4(^tS\Omega S)(\Z)$, therefore $\Gamma(f,g)\subset\Sp_4(\Omega)(\Z)$ is arithmetic if and only if $\Gamma(\underbar{f}, \underbar{g})\subset\Sp_4(^tS\Omega S)(\Z)$ is arithmetic.
\end{rmk}

In Table \ref{table:newarithmetic}, we list the $15$ pairs $f,g$ of \cite[Table 2]{SV} which we show yield arithmetic monodromy groups in Section \ref{proofarithmetic}  (cf. Remark \ref{criterion}). In fact, we obtain the following theorem:

\begin{theorem}\label{maintheorem}
 The monodromy groups associated to the $15$ pairs of polynomials in Table \ref{table:newarithmetic} are arithmetic.
\end{theorem}

Finally, in Table \ref{table:unknown}, we list the remaining $11$ pairs $f,g$ of \cite[Table 2]{SV}, for which, the question to determine the arithmeticity or thinness of the associated monodromy groups is still open. By using Remark \ref{remarkf(-x)g(-x)}, we observe that it is enough to prove the arithmeticity or thinness for Examples \ref{unknown1}-\ref{unknown6} of Table \ref{table:unknown}, and the same will follow for the rest of the examples.

We now note the following remarks:
\begin{rmk}
 Brav and Thomas \cite{BT} prove the thinness of $7$ monodromy groups associated to Examples \ref{thin1}-\ref{thin7} of Table \ref{table:thin}. They use a ping-pong argument to show that the monodromy groups associated to these pairs are either a free group or contain a free subgroup of finite index; and by comparing the cohomological dimensions of the respective subgroups with that of $\Sp_4(\Omega)(\Z)$, they show that these groups can not be of finite index. We hope that the same method could be applied to prove thinness for some of the examples in Table \ref{table:unknown}. 
\end{rmk}

\begin{rmk}\label{criterion} To prove Theorem \ref{maintheorem}, it is enough to prove the arithmeticity of $8$ monodromy groups associated to Examples \ref{arithmetic18}-\ref{arithmetic20} of Table \ref{table:newarithmetic}, by using Remark \ref{remarkf(-x)g(-x)}. The method to prove the arithmeticity for the examples of Table \ref{table:newarithmetic} is same as that of \cite{SS} and \cite{SV}. In particular, we show that the monodromy groups $\Gamma(f,g)$ of Table \ref{table:newarithmetic} intersect the $\Z$ points  $U(\Z)$ of the unipotent radical $U$ of a Borel subgroup $B$ of $\Sp_4(\Omega)$ in a finite index subgroup of $U(\Z)$. Since the monodromy groups of Table \ref{table:newarithmetic} are Zariski dense in $\Sp_4(\Omega)$ by \cite{BH}, their arithmeticity follows by \cite{T}. 
 \end{rmk}

 \begin{rmk}\label{Venkataramana}
 It follows from \cite[Theorem 3.5]{Ve} that if $\Gamma$ is a Zariski dense subgroup of $\Sp_4(\Omega)(\Z)$, and intersects the highest  and second highest root groups non-trivially, then  $\Gamma$ has finite index in $\Sp_4(\Omega)(\Z)$. Note that, once we show that the groups $\Gamma(f,g)$ associated to the pairs in Table \ref{table:newarithmetic}, intersect the subgroup $U(\Z)$ of $\Sp_4(\Omega)(\Z)$ in a finite index subgroup of $U(\Z)$, it  follows automatically that $\Gamma(f,g)$ also intersects the highest  and second highest root groups non-trivially, and the arithmeticity of these groups also follows from \cite[Theorem 3.5]{Ve}.
\end{rmk}

\section{Tables}
In this section, we split \cite[Table 2]{SV} in $4$ subtables, that is, in Tables \ref{table:arithmetic}, \ref{table:thin}, \ref{table:newarithmetic}, and \ref{table:unknown}, depending on the progress to answer the question to determine the pairs of polynomials in \cite[Table 2]{SV}, for which, the associated monodromy groups are arithmetic or thin. 

In Table \ref{table:arithmetic}, we list the $12$ {\it arithmetic} monodromy groups of \cite[Table 2]{SV}, for which, the arithmeticity follows from \cite{SS} and \cite{SV} (cf. Remark \ref{remarkf(-x)g(-x)}). In Table \ref{table:thin}, we list the $13$ {\it thin} monodromy groups of \cite[Table 2]{SV}, for which, the thinness follows from \cite{BT} (cf. Remark \ref{remarkf(-x)g(-x)}). 
In Table \ref{table:newarithmetic}, we list the $15$ {\it arithmetic} monodromy groups of \cite[Table 2]{SV}, for which, the arithmeticity is shown in Section \ref{proofarithmetic} of this article (cf. Remark \ref{remarkf(-x)g(-x)}).
Finally, in Table \ref{table:unknown}, we list the remaining $11$ monodromy groups of \cite[Table 2]{SV}, for which, the {\it arithmeticity} or {\it thinness} is {\it unknown}. In fact, it is enough to prove the arithmeticity or thinness of Examples \ref{unknown1}-\ref{unknown6} in Table \ref{table:unknown}, and the same will follow for the rest of the examples (cf. Remark \ref{remarkf(-x)g(-x)}). 

In Tables \ref{table:arithmetic}, \ref{table:thin}, \ref{table:newarithmetic}, and \ref{table:unknown}, we identify the pair $(f(-X),g(-X))$ (up to a transposition) for a given pair $(f(X),g(X))$, by putting the number of the earlier pair inside a bracket following the number of the later one; and when the pairs $(f(X),g(X))$ and $(f(-X),g(-X))$ are same (up to a transposition), we put the same number of the pair in a following bracket.

{\renewcommand{\arraystretch}{1.5}
\begin{table}[h]
\tiny
\addtolength{\tabcolsep}{-4pt}
\caption{List of the $12$ {\it arithmetic} monodromy groups of \cite[Table 2]{SV}, for which, the arithmeticity follows from \cite{SS} and \cite{SV} (cf. Remark \ref{remarkf(-x)g(-x)}).}
\newcounter{rownum-2}
\setcounter{rownum-2}{0}
\centering
\begin{tabular}{ |c|  c|   c| c| c| c|}
\hline

  No. & $f(X)$ & $g(X)$ & $\alpha$ & $\beta$ & $f(X)-g(X)$ \\ \hline
  \refstepcounter{rownum-2}\arabic{rownum-2}\label{arithmetic1} & $X^4-4X^3+6X^2-4X+1$ &$X^4+2X^3+3X^2+2X+1$  &$0,0,0,0$ &$\frac{1}{3},\frac{1}{3},\frac{2}{3},\frac{2}{3}$ &$-6X^3+3X^2-6X$ \\ \hline
  
   \refstepcounter{rownum-2}\arabic{rownum-2}\label{arithmetic2} & $X^4-4X^3+6X^2-4X+1$ &$X^4+2X^2+1$  &$0,0,0,0$ &$\frac{1}{4},\frac{1}{4}$,$\frac{3}{4}$,$\frac{3}{4}$ &$-4X^3+4X^2-4X$ \\ \hline
   
   \refstepcounter{rownum-2}\arabic{rownum-2}\label{arithmetic3} & $X^4-4X^3+6X^2-4X+1$ &$X^4+X^3+2X^2+X+1$  &$0,0,0,0$ &$\frac{1}{3}$,$\frac{2}{3}$,$\frac{1}{4}$,$\frac{3}{4}$ &$-5X^3+4X^2-5X$ \\ \hline
   
   \refstepcounter{rownum-2}\arabic{rownum-2}\label{arithmetic4} & $X^4-4X^3+6X^2-4X+1$ &$X^4+X^2+1$  &$0,0,0,0$ &$\frac{1}{3}$,$\frac{2}{3}$,$\frac{1}{6}$,$\frac{5}{6}$ &$-4X^3+5X^2-4X$ \\ \hline
   
   \refstepcounter{rownum-2}\arabic{rownum-2}\label{arithmetic5} & $ X^4-4X^3+6X^2-4X+1$ &$ X^4-X^3+2X^2-X+1$  &$ 0,0,0,0$ &$\frac{1}{4},\frac{3}{4},\frac{1}{6},\frac{5}{6}$ &$-3X^3+4X^2-3X$ \\ \hline
   
   \refstepcounter{rownum-2}\arabic{rownum-2}\label{arithmetic6} & $ X^4-4X^3+6X^2-4X+1$ &$ X^4-X^3+X^2-X+1$  &$ 0,0,0,0$ &$\frac{1}{10},\frac{3}{10},\frac{7}{10},\frac{9}{10}$ &$-3X^3+5X^2-3X$ \\ \hline
   
   \refstepcounter{rownum-2}\arabic{rownum-2}\label{arithmetic10}(\ref{arithmetic1}) & $X^4+4X^3+6X^2+4X+1$ &$X^4-2X^3+3X^2-2X+1$  &$\frac{1}{2}$,$\frac{1}{2}$,$\frac{1}{2}$,$\frac{1}{2}$ &$\frac{1}{6}$,$\frac{1}{6}$,$\frac{5}{6}$,$\frac{5}{6}$ &$6X^3+3X^2+6X$ \\ \hline
   
    \refstepcounter{rownum-2}\arabic{rownum-2}\label{arithmetic7}(\ref{arithmetic2}) & $X^4+4X^3+6X^2+4X+1$ &$X^4+2X^2+1$  &$\frac{1}{2}$,$\frac{1}{2}$,$\frac{1}{2}$,$\frac{1}{2}$ &$\frac{1}{4}$,$\frac{1}{4}$,$\frac{3}{4}$,$\frac{3}{4}$ &$4X^3+4X^2+4X$ \\ \hline
    
  \refstepcounter{rownum-2}\arabic{rownum-2}\label{arithmetic12}(\ref{arithmetic3}) & $X^4+4X^3+6X^2+4X+1$ &$X^4-X^3+2X^2-X+1$  &$\frac{1}{2}$,$\frac{1}{2}$,$\frac{1}{2}$,$\frac{1}{2}$ &$\frac{1}{4}$,$\frac{3}{4}$,$\frac{1}{6}$,$\frac{5}{6}$ &$5X^3+4X^2+5X$ \\ \hline
    
  \refstepcounter{rownum-2}\arabic{rownum-2}\label{arithmetic11}(\ref{arithmetic4}) & $X^4+4X^3+6X^2+4X+1$ &$X^4+X^2+1$  &$\frac{1}{2}$,$\frac{1}{2}$,$\frac{1}{2}$,$\frac{1}{2}$ &$\frac{1}{3}$,$\frac{2}{3}$,$\frac{1}{6}$,$\frac{5}{6}$ &$4X^3+5X^2+4X$ \\ \hline
  
  \refstepcounter{rownum-2}\arabic{rownum-2}\label{arithmetic8}(\ref{arithmetic5}) & $X^4+4X^3+6X^2+4X+1$ &$X^4+X^3+2X^2+X+1$  &$\frac{1}{2}$,$\frac{1}{2}$,$\frac{1}{2}$,$\frac{1}{2}$ &$\frac{1}{4}$,$\frac{3}{4}$,$\frac{1}{3}$,$\frac{2}{3}$ &$3X^3+4X^2+3X$ \\ \hline
  
  \refstepcounter{rownum-2}\arabic{rownum-2}\label{arithmetic9}(\ref{arithmetic6}) & $X^4+4X^3+6X^2+4X+1$ &$X^4+X^3+X^2+X+1$  &$\frac{1}{2}$,$\frac{1}{2}$,$\frac{1}{2}$,$\frac{1}{2}$ &$\frac{1}{5}$,$\frac{2}{5}$,$\frac{3}{5}$,$\frac{4}{5}$ &$3X^3+5X^2+3X$ \\ \hline
  \end{tabular}
\label{table:arithmetic}
\end{table}}
\vspace{1cm}
{\renewcommand{\arraystretch}{1.5}   
{\tiny\begin{table}[h]
\addtolength{\tabcolsep}{-4pt}
\caption{List of the $13$ {\it thin} monodromy groups of \cite[Table 2]{SV}, for which, the thinness follows from \cite{BT} (cf. Remark \ref{remarkf(-x)g(-x)}).}
\newcounter{rownum-4}
\setcounter{rownum-4}{0}
\centering
\begin{tabular}{ |c|  c|   c| c| c| c|}
\hline
 No. & $f(X)$ & $g(X)$ & $\alpha$ & $\beta$ & $f(X)-g(X)$\\ \hline
 
  \refstepcounter{rownum-4}\arabic{rownum-4}\label{thin1}(\ref{thin1}) & $X^4-4X^3+6X^2-4X+1$ &$ X^4+4X^3+6X^2+4X+1$  &$0,0,0,0$ &$\frac{1}{2}, \frac{1}{2},\frac{1}{2},\frac{1}{2}$ &$-8X^3-8X$ \\ \hline

  \refstepcounter{rownum-4}\arabic{rownum-4}\label{thin2} &  $X^4-4X^3+6X^2-4X+1$ &$ X^4+3X^3+4X^2+3X+1$  & $0,0,0,0$ &$\frac{1}{2}, \frac{1}{2},\frac{1}{3},\frac{2}{3}$ &$ -7X^3+2X^2-7X$ \\ \hline
  
  \refstepcounter{rownum-4}\arabic{rownum-4}\label{thin3} & $ X^4-4X^3+6X^2-4X+1$ &$X^4+2X^3+2X^2+2X+1$  &$0,0,0,0$ &$\frac{1}{2}, \frac{1}{2},\frac{1}{4},\frac{3}{4}$ &$ -6X^3+4X^2-6X$ \\ \hline
  
  \refstepcounter{rownum-4}\arabic{rownum-4}\label{thin4} & $ X^4-4X^3+6X^2-4X+1$ & $ X^4+X^3+X^2+X+1$  &$0,0,0,0$ &$\frac{1}{5},\frac{2}{5},\frac{3}{5},\frac{4}{5}$ &$-5X^3+5X^2-5X$ \\ \hline
  
  \refstepcounter{rownum-4}\arabic{rownum-4}\label{thin5} & $ X^4-4X^3+6X^2-4X+1$ &$ X^4+X^3+X+1$  &$0,0,0,0$ &$\frac{1}{2},\frac{1}{2},\frac{1}{6},\frac{5}{6}$ &$-5X^3+6X^2-5X$ \\ \hline
  
  \refstepcounter{rownum-4}\arabic{rownum-4}\label{thin6} & $X^4-4X^3+6X^2-4X+1$ &$ X^4+1$  & $0,0,0,0 $&$\frac{1}{8},\frac{3}{8},\frac{5}{8},\frac{7}{8}$ &$-4X^3+6X^2-4X$ \\ \hline
  
  \refstepcounter{rownum-4}\arabic{rownum-4}\label{thin7} & $ X^4-4X^3+6X^2-4X+1$ &$ X^4-X^2+1$  &$0,0,0,0$ &$\frac{1}{12},\frac{5}{12},\frac{7}{12},\frac{11}{12}$ &$ -4X^3+7X^2-4X$ \\ \hline
  
  \refstepcounter{rownum-4}\arabic{rownum-4}\label{thin10}(\ref{thin2}) & $X^4+4X^3+6X^2+4X+1$ &$X^4-3X^3+4X^2-3X+1$  &$\frac{1}{2}$,$\frac{1}{2}$,$\frac{1}{2}$,$\frac{1}{2}$ &0,0,$\frac{1}{6}$,$\frac{5}{6}$ &$7X^3+2X^2+7X$ \\ \hline
  
   \refstepcounter{rownum-4}\arabic{rownum-4}\label{thin9}(\ref{thin3}) & $X^4+4X^3+6X^2+4X+1$ &$X^4-2X^3+2X^2-2X+1$  &$\frac{1}{2}$,$\frac{1}{2}$,$\frac{1}{2}$,$\frac{1}{2}$ &$0,0,\frac{1}{4}$,$\frac{3}{4}$ &$6X^3+4X^2+6X$ \\ \hline
   
    \refstepcounter{rownum-4}\arabic{rownum-4}\label{thin12}(\ref{thin4}) & $X^4+4X^3+6X^2+4X+1$ &$X^4-X^3+X^2-X+1$  &$\frac{1}{2}$,$\frac{1}{2}$,$\frac{1}{2}$,$\frac{1}{2}$ &$\frac{1}{10}$,$\frac{3}{10}$,$\frac{7}{10}$,$\frac{9}{10}$ &$5X^3+5X^2+5X$ \\ \hline
      
   \refstepcounter{rownum-4}\arabic{rownum-4}\label{thin8}(\ref{thin5}) & $X^4+4X^3+6X^2+4X+1$ &$X^4-X^3-X+1$  &$\frac{1}{2}$,$\frac{1}{2}$,$\frac{1}{2}$,$\frac{1}{2}$ &0,0,$\frac{1}{3}$,$\frac{2}{3}$ &$5X^3+6X^2+5X$ \\ \hline
    
   \refstepcounter{rownum-4}\arabic{rownum-4}\label{thin11}(\ref{thin6}) & $X^4+4X^3+6X^2+4X+1$ &$X^4+1$  &$\frac{1}{2}$,$\frac{1}{2}$,$\frac{1}{2}$,$\frac{1}{2}$ &$\frac{1}{8}$,$\frac{3}{8}$,$\frac{5}{8}$,$\frac{7}{8}$ &$4X^3+6X^2+4X$ \\ \hline
   
  \refstepcounter{rownum-4}\arabic{rownum-4}\label{thin13}(\ref{thin7}) & $X^4+4X^3+6X^2+4X+1$ &$X^4-X^2+1$  &$\frac{1}{2}$,$\frac{1}{2}$,$\frac{1}{2}$,$\frac{1}{2}$ &$\frac{1}{12}$,$\frac{5}{12}$,$\frac{7}{12}$,$\frac{11}{12}$ &$4X^3+7X^2+4X$ \\ \hline
  
\end{tabular}\label{table:thin}
\end{table}}}
 
\newpage
{\renewcommand{\arraystretch}{1.37}
\begin{table}[h]
\tiny
\addtolength{\tabcolsep}{-4pt}
\caption{List of the $15$ {\it arithmetic} monodromy groups of \cite[Table 2]{SV}, for which, the arithmeticity is shown in Section \ref{proofarithmetic} (cf. Remark \ref{remarkf(-x)g(-x)}).}
\newcounter{rownum-3}
\setcounter{rownum-3}{0}
\centering
\begin{tabular}{ |c|  c|   c| c| c| c|}
\hline

  No. & $f(X)$ & $g(X)$ & $\alpha$ & $\beta$ & $f(X)-g(X)$ \\ \hline
   
    \refstepcounter{rownum-3}\arabic{rownum-3}\label{arithmetic18} & $X^4+3X^3+4X^2+3X+1$ &$X^4+2X^2+1$  &$\frac{1}{2}$,$\frac{1}{2}$,$\frac{1}{3}$,$\frac{2}{3}$  &$\frac{1}{4}$,$\frac{1}{4}$,$\frac{3}{4}$,$\frac{3}{4}$ &$3X^3+2X^2+3X$ \\ \hline
    
\refstepcounter{rownum-3}\arabic{rownum-3}\label{arithmetic13} & $X^4+2X^3+3X^2+2X+1$ &$X^4-2X^3+2X^2-2X+1$  &$\frac{1}{3}$,$\frac{1}{3}$,$\frac{2}{3}$,$\frac{2}{3}$ &0,0,$\frac{1}{4}$,$\frac{3}{4}$ &$4X^3+X^2+4X$ \\ \hline
    
     \refstepcounter{rownum-3}\arabic{rownum-3}\label{arithmetic14} & $X^4+2X^3+3X^2+2X+1$ &$X^4-3X^3+4X^2-3X+1$  &$\frac{1}{3}$,$\frac{1}{3}$,$\frac{2}{3}$,$\frac{2}{3}$ &0,0,$\frac{1}{6}$,$\frac{5}{6}$ &$5X^3-X^2+5X$ \\ \hline
  
  \refstepcounter{rownum-3}\arabic{rownum-3}\label{arithmetic15}(\ref{arithmetic15}) & $X^4+2X^3+3X^2+2X+1$ &$X^4-2X^3+3X^2-2X+1$  &$\frac{1}{3}$,$\frac{1}{3}$,$\frac{2}{3}$,$\frac{2}{3}$ &$\frac{1}{6}$,$\frac{1}{6}$,$\frac{5}{6}$,$\frac{5}{6}$ &$4X^3+4X$ \\ \hline
  
  \refstepcounter{rownum-3}\arabic{rownum-3}\label{arithmetic16} & $X^4+2X^3+3X^2+2X+1$ &$X^4-X^3+2X^2-X+1$  &$\frac{1}{3}$,$\frac{1}{3}$,$\frac{2}{3}$,$\frac{2}{3}$ &$\frac{1}{4}$,$\frac{3}{4}$,$\frac{1}{6}$,$\frac{5}{6}$ &$3X^3+X^2+3X$ \\ \hline
  
    \refstepcounter{rownum-3}\arabic{rownum-3}\label{arithmetic17} & $X^4+2X^3+3X^2+2X+1$ &$X^4-X^3+X^2-X+1$  &$\frac{1}{3}$,$\frac{1}{3}$,$\frac{2}{3}$,$\frac{2}{3}$ &$\frac{1}{10}$,$\frac{3}{10}$,$\frac{7}{10}$,$\frac{9}{10}$ &$3X^3+2X^2+3X$ \\ \hline
       
   \refstepcounter{rownum-3}\arabic{rownum-3}\label{arithmetic19} & $X^4+3X^3+4X^2+3X+1$ &$X^4-X^3+2X^2-X+1$  &$\frac{1}{2}$,$\frac{1}{2}$,$\frac{1}{3}$,$\frac{2}{3}$  &$\frac{1}{4}$,$\frac{3}{4}$,$\frac{1}{6}$,$\frac{5}{6}$ &$4X^3+2X^2+4X$ \\ \hline
    
     \refstepcounter{rownum-3}\arabic{rownum-3}\label{arithmetic20} & $X^4+3X^3+4X^2+3X+1$ &$X^4-X^2+1$  &$\frac{1}{2}$,$\frac{1}{2}$,$\frac{1}{3}$,$\frac{2}{3}$  &$\frac{1}{12}$,$\frac{5}{12}$,$\frac{7}{12}$,$\frac{11}{12}$ &$3X^3+5X^2+3X$ \\ \hline
     
         \refstepcounter{rownum-3}\arabic{rownum-3}\label{arithmetic22}(\ref{arithmetic18}) & $X^4+2X^2+1$ &$X^4-3X^3+4X^2-3X+1$  &$\frac{1}{4}$,$\frac{1}{4}$,$\frac{3}{4}$,$\frac{3}{4}$  &0,0,$\frac{1}{6}$,$\frac{5}{6}$ &$3X^3-2X^2+3X$ \\ \hline
     
      \refstepcounter{rownum-3}\arabic{rownum-3}\label{arithmetic23}(\ref{arithmetic13}) & $X^4+2X^3+2X^2+2X+1$ &$X^4-2X^3+3X^2-2X+1$  &$\frac{1}{2}$,$\frac{1}{2}$,$\frac{1}{4}$,$\frac{3}{4}$  &$\frac{1}{6}$,$\frac{1}{6}$,$\frac{5}{6}$,$\frac{5}{6}$ &$4X^3-X^2+4X$ \\ \hline
      
       \refstepcounter{rownum-3}\arabic{rownum-3}\label{arithmetic21}(\ref{arithmetic14}) & $X^4+3X^3+4X^2+3X+1$ &$X^4-2X^3+3X^2-2X+1$  &$\frac{1}{2}$,$\frac{1}{2}$,$\frac{1}{3}$,$\frac{2}{3}$  &$\frac{1}{6}$,$\frac{1}{6}$,$\frac{5}{6}$,$\frac{5}{6}$ &$5X^3+X^2+5X$ \\ \hline
           
\refstepcounter{rownum-3}\arabic{rownum-3}\label{arithmetic25}(\ref{arithmetic16}) & $ X^4+X^3+2X^2+X+1$ &$ X^4-2X^3+3X^2-2X+1$  &$\frac 13,\frac 23,\frac 14,\frac 34$ & $\frac 16,\frac 16,\frac 56,\frac 56$ &$3X^3-X^2+3X$ \\ \hline

\refstepcounter{rownum-3}\arabic{rownum-3}\label{arithmetic26}(\ref{arithmetic17}) & $X^4+X^3+X^2+X+1$ &$X^4-2X^3+3X^2-2X+1$  &$\frac{1}{5}$,$\frac{2}{5}$,$\frac{3}{5}$,$\frac{4}{5}$  &$\frac{1}{6}$,$\frac{1}{6}$,$\frac{5}{6}$,$\frac{5}{6}$ &$3X^3-2X^2+3X$ \\ \hline

      \refstepcounter{rownum-3}\arabic{rownum-3}\label{arithmetic24}(\ref{arithmetic19}) & $X^4+X^3+2X^2+X+1$ &$X^4-3X^3+4X^2-3X+1$  &$\frac{1}{3}$,$\frac{2}{3}$,$\frac{1}{4}$,$\frac{3}{4}$  &0,0,$\frac{1}{6}$,$\frac{5}{6}$ &$4X^3-2X^2+4X$ \\ \hline
        
\refstepcounter{rownum-3}\arabic{rownum-3}\label{arithmetic27}(\ref{arithmetic20}) & $X^4-3X^3+4X^2-3X+1$ &$X^4-X^2+1$  &0,0,$\frac{1}{6}$,$\frac{5}{6}$  &$\frac{1}{12}$,$\frac{5}{12}$,$\frac{7}{12}$,$\frac{11}{12}$ &$-3X^3+5X^2-3X$ \\ \hline
\end{tabular}
\label{table:newarithmetic}
\end{table}}
\vspace{1cm}
{\renewcommand{\arraystretch}{1.37}   
{\tiny\begin{table}[h]
\addtolength{\tabcolsep}{-4pt}
\caption{List of the remaining $11$ monodromy groups of \cite[Table 2]{SV}, for which, the {\it arithmeticity} or {\it thinness} is {\it unknown}.}
\newcounter{rownum-5}
\setcounter{rownum-5}{0}
\centering
\begin{tabular}{ |c|  c|   c| c| c| c|}
\hline
 No. & $f(X)$ & $g(X)$ & $\alpha$ & $\beta$ & $f(X)-g(X)$\\ \hline
 
   \refstepcounter{rownum-5}\arabic{rownum-5}\label{unknown1} & $X^4-X^3-X+1$ &$X^4+2X^3+2X^2+2X+1$  &0,0,$\frac{1}{3}$,$\frac{2}{3}$ &$\frac{1}{2}$,$\frac{1}{2}$,$\frac{1}{4}$,$\frac{3}{4}$ &$-3X^3-2X^2-3X$ \\ \hline
   
     \refstepcounter{rownum-5}\arabic{rownum-5}\label{unknown2} & $X^4+3X^3+4X^2+3X+1$ &$X^4-2X^3+2X^2-2X+1$  &$\frac{1}{2}$,$\frac{1}{2}$,$\frac{1}{3}$,$\frac{2}{3}$ &0,0,$\frac{1}{4}$,$\frac{3}{4}$ &$5X^3+2X^2+5X$ \\ \hline
     
     \refstepcounter{rownum-5}\arabic{rownum-5}\label{unknown3}(\ref{unknown3}) & $X^4+3X^3+4X^2+3X+1$ &$X^4-3X^3+4X^2-3X+1$  &$\frac{1}{2}$,$\frac{1}{2}$,$\frac{1}{3}$,$\frac{2}{3}$  &0,0,$\frac{1}{6}$,$\frac{5}{6}$ &$6X^3+6X$ \\ \hline
     
     \refstepcounter{rownum-5}\arabic{rownum-5}\label{unknown4} & $X^4+3X^3+4X^2+3X+1$ &$X^4+1$  &$\frac{1}{2}$,$\frac{1}{2}$,$\frac{1}{3}$,$\frac{2}{3}$  &$\frac{1}{8}$,$\frac{3}{8}$,$\frac{5}{8}$,$\frac{7}{8}$ &$3X^3+4X^2+3X$ \\ \hline
  
  \refstepcounter{rownum-5}\arabic{rownum-5}\label{unknown5} & $X^4+3X^3+4X^2+3X+1$ &$X^4-X^3+X^2-X+1$  &$\frac{1}{2}$,$\frac{1}{2}$,$\frac{1}{3}$,$\frac{2}{3}$  &$\frac{1}{10}$,$\frac{3}{10}$,$\frac{7}{10}$,$\frac{9}{10}$ &$4X^3+3X^2+4X$ \\ \hline
  
   \refstepcounter{rownum-5}\arabic{rownum-5}\label{unknown6} & $X^4-2X^3+2X^2-2X+1$ &$X^4+X^3+X^2+X+1$  &0,0,$\frac{1}{4}$,$\frac{3}{4}$ &$\frac{1}{5}$,$\frac{2}{5}$,$\frac{3}{5}$,$\frac{4}{5}$ &$-3X^3+X^2-3X$ \\ \hline
  
  \refstepcounter{rownum-5}\arabic{rownum-5}\label{unknown7}(\ref{unknown1}) & $X^4-2X^3+2X^2-2X+1$ &$X^4+X^3+X+1$  &0,0,$\frac{1}{4}$,$\frac{3}{4}$ &$\frac{1}{2}$,$\frac{1}{2}$,$\frac{1}{6}$,$\frac{5}{6}$ &$-3X^3+2X^2-3X$ \\ \hline
  
   \refstepcounter{rownum-5}\arabic{rownum-5}\label{unknown8}(\ref{unknown2}) & $X^4+2X^3+2X^2+2X+1$ &$X^4-3X^3+4X^2-3X+1$  &$\frac{1}{2}$,$\frac{1}{2}$,$\frac{1}{4}$,$\frac{3}{4}$  &0,0,$\frac{1}{6}$,$\frac{5}{6}$ &$5X^3-2X^2+5X$ \\ \hline
      
   \refstepcounter{rownum-5}\arabic{rownum-5}\label{unknown9}(\ref{unknown4}) & $X^4-3X^3+4X^2-3X+1$ &$X^4+1$  &0,0,$\frac{1}{6}$,$\frac{5}{6}$  &$\frac{1}{8}$,$\frac{3}{8}$,$\frac{5}{8}$,$\frac{7}{8}$ &$-3X^3+4X^2-3X$ \\ \hline
   
   \refstepcounter{rownum-5}\arabic{rownum-5}\label{unknown10}(\ref{unknown5}) & $X^4+X^3+X^2+X+1$ &$X^4-3X^3+4X^2-3X+1$  &$\frac{1}{5}$,$\frac{2}{5}$,$\frac{3}{5}$,$\frac{4}{5}$  &0,0,$\frac{1}{6}$,$\frac{5}{6}$ &$4X^3-3X^2+4X$ \\ \hline
   
   \refstepcounter{rownum-5}\arabic{rownum-5}\label{unknown11}(\ref{unknown6}) & $X^4+2X^3+2X^2+2X+1$ &$X^4-X^3+X^2-X+1$  &$\frac{1}{2}$,$\frac{1}{2}$,$\frac{1}{4}$,$\frac{3}{4}$  &$\frac{1}{10}$,$\frac{3}{10}$,$\frac{7}{10}$,$\frac{9}{10}$ &$3X^3+X^2+3X$ \\ \hline
   
   \end{tabular}\label{table:unknown}
\end{table}}}

\section{Proof of Theorem \ref{maintheorem}}\label{proofarithmetic}

We will first compute the symplectic form $\Omega$ (up to scalar multiples) preserved by each monodromy group in Table \ref{table:newarithmetic}, then show that there exists a basis $\{\e_1,\e_2,\e_2^*,\e_1^*\}$ of $\Q^4$, with respect to which, the matrix form of $\Omega$ is anti-diagonal. It can be checked easily that, with respect to the basis $\{\e_1,\e_2,\e_2^*,\e_1^*\}$, the diagonal matrices in $\Sp_4(\Omega)$ form a maximal torus, the group of upper triangular matrices in $\Sp_4(\Omega)$ form a Borel subgroup $\mathrm{B}$, and the group of unipotent upper triangular matrices in $\Sp_4(\Omega)$ form the unipotent radical $\mathrm{U}$ of $\mathrm{B}$.

Note that $\mathrm{U}$ is a nilpotent subgroup of $\GL_4(\R)$, and it follows from \cite[Theorem 2.1]{Rag} that if $\Gamma(f,g)\cap\mathrm{U}(\Z)$ is a Zariski dense subgroup of $\mathrm{U}$ then $\mathrm{U}/\Gamma(f,g)\cap\mathrm{U}(\Z)$ is compact, and hence $\Gamma(f,g)\cap\mathrm{U}(\Z)$ has finite index in $\mathrm{U}(\Z)$.  Therefore, to show that the monodromy group $\Gamma(f,g)$, with respect to the basis $\{\e_1,\e_2,\e_2^*,\e_1^*\}$, intersects $\mathrm{U}(\Z)$ in a finite index subgroup of $\mathrm{U}(\Z)$ (cf. Remark \ref{criterion}), it is enough to show that $\Gamma(f,g)\cap\mathrm{U}(\Z)$ is Zariski dense in $\mathrm{U}$; and to show this, it is enough to show that $\Gamma(f,g)$ contains non-trivial unipotent elements corresponding to each of the positive roots; and the proof  of arithmeticity of $\Gamma(f,g)$ follows from \cite{T} (cf. \cite[Theorem 3.5]{Ve}, Remark \ref{Venkataramana}).

We now explain the logic behind finding the unipotent elements in $\Gamma(f,g)$. For that, let $\{e_1,e_2,e_3,e_4\}$ be the standard basis of $\Q^4$ over $\Q$, and $v$ be the last column vector of the matrix $\left(A^{-1}B-I\right)$, that is, $v=\left(A^{-1}B-I\right)(e_4)$, where $I$ is the $4\times4$ identity matrix. It also has been remarked in \cite[Remark 5.1]{SV} that if we could find an element $\gamma\in\Gamma(f,g)$, for which, the absolute value of the coefficient of $e_4$ in $\gamma(v)$, is non-zero and $\leq2$, then the method of the proof of \cite[Theorem 1.1]{SV} can be applied, and the arithmeticity of $\Gamma(f,g)$ will follow. 

Actually, in the proof of \cite[Theorem 1.1]{SV} (for the case $n=4$), the element $\gamma$ is $A^k$ for some $k\geq1$, for which,  the absolute value of the coefficient of $e_4$ in $A^k(v)$, is non-zero and $\leq2$, and it has been shown that the subgroup $\Gamma_r\subset\Gamma(f,g)$ generated by $C=A^{-1}B$, $A^{-k}CA^{k}$ and $A^kCA^{-k}$, contains all the required unipotent elements. Once we find the element $\gamma\in\Gamma(f,g)$ which satisfy the above condition, replace $A^k$ by $\gamma$, and exactly the same method of the proof of \cite[Theorem 1.1]{SV} is applied to compute the unipotent elements. We did some experiments, and found such $\gamma$ in $7$ examples of Table \ref{table:newarithmetic}, and applied the method of the proof of \cite[Theorem 1.1]{SV}.  

There are many cases where such $\gamma\in\Gamma(f,g)$ can not even exist: if some $\Gamma(f,g)$ is going to be thin, it is not possible to find the element $\gamma$ (because the existence of such $\gamma\in\Gamma(f,g)$ will show the arithmeticity of $\Gamma(f,g)$); and if the greatest common divisor of the coefficients of $e_1,e_2,e_3,e_4$ in $v=\left(A^{-1}B-I\right)(e_4)$ is $\geq3$ (cf. Examples \ref{arithmetic1} and \ref{arithmetic2} of Table \ref{table:arithmetic}, and Example \ref{arithmetic15} of Table \ref{table:newarithmetic}), then also it is not possible to find the element $\gamma$ (because the coefficients of $e_1,e_2,e_3,e_4$ in $\gamma(v)$ will always be multiple of an integer $\geq3$).

An interesting example is Example \ref{arithmetic15} of Table \ref{table:newarithmetic}, where the greatest common divisor of the coefficients of $e_1,e_2,e_3,e_4$ in $v=\left(A^{-1}B-I\right)(e_4)$ is $4$, and therefore such $\gamma\in\Gamma(f,g)$ can not exist but we are still able to prove the arithmeticity of the associated monodromy group $\Gamma(f,g)$ (cf. Subsection \ref{proofarithmetic15}). 

We now proceed to prove the arithmeticity of the monodromy groups associated to the examples in Table \ref{table:newarithmetic}.

\subsection{Example \ref{arithmetic18} of Table \ref{table:newarithmetic}}\label{proofarithmetic18}
This is Example 33 of \cite[Table 2]{SV}. In this case  $$\alpha=\left(\frac{1}{2},\frac{1}{2},\frac{1}{3},\frac{2}{3}\right), \quad \beta=\left(\frac{1}{4},\frac{1}{4},\frac{3}{4},\frac{3}{4}\right);$$  $$f(X)=X^4+3X^3+4X^2+3X+1, \quad g(X)=X^4+2X^2+1;$$ and $f(X)-g(X)=3X^3+2X^2+3X$.

Let $A$ and $B$  be the companion  matrices of $f(X)$ and  $g(X)$ respectively, and let $C=A^{-1}B$. Then
{\scriptsize \[A=\begin{pmatrix}  \begin {array}{rrrr} 0&0&0&-1\\ \noalign{\medskip}1&0&0&-3
\\ \noalign{\medskip}0&1&0&-4\\ \noalign{\medskip}0&0&1&-3\end {array}
  \end{pmatrix},\quad  B=\begin{pmatrix} \begin {array}{cccr} 0&0&0&-1\\ \noalign{\medskip}1&0&0&0
\\ \noalign{\medskip}0&1&0&-2\\ \noalign{\medskip}0&0&1&0\end {array}
 \end{pmatrix},\quad C=\begin{pmatrix} \begin {array}{rrrr} 1&0&0&3\\ \noalign{\medskip}0&1&0&2
\\ \noalign{\medskip}0&0&1&3\\ \noalign{\medskip}0&0&0&1\end {array}
 \end{pmatrix}.\]}   Let $\Gamma(f,g)$ be the  subgroup  of
$\SL_4(\Z)$ generated by $A$ and $B$. 
\subsection*{The invariant symplectic form} Let us denote $\Omega(v_1,v_2)$ by $v_1.v_2$, for any pairs of vectors $v_1,v_2\in\Q^4$. Let $\{e_1,e_2,e_3,e_4\}$ be the standard basis of $\Q^4$ over $\Q$, and $v=3e_1+2e_2+3e_3$, which is the last column vector of $C-\mathrm{I}$, where $\mathrm{I}$ is the $4\times 4$ identity matrix. Since $C$ preserves the form $\Omega$, for $1\leq i\leq 3$, we obtain
\begin{align*}
e_i.e_4&=e_i.(3e_1+2e_2+3e_3+e_4)\\
&=e_i.(v+e_4)\\
&=e_i.v+e_i.e_4.
\end{align*}
This implies that 
\begin{equation}\label{orthogonalv}
 e_i.v=0\qquad \mbox{for}\ 1\leq i\leq 3.
\end{equation}
That is, $v$ is $\Omega$- orthogonal to the vectors $e_1, e_2, e_3$ and $e_4.v\neq 0$ (since $\Omega$ is non-degenerate). Since $B$ preserves $\Omega$, we obtain
\begin{align*}
e_1.e_2=e_2.e_3=e_3.e_4&=e_4.(-e_1-2e_3)\\
&=e_1.e_4+2e_3.e_4.
\end{align*}
This implies that
\begin{equation}\label{equation2}
e_1.e_4=-e_3.e_4.
\end{equation}
It now follows from (\ref{orthogonalv}) and the invariance of $\Omega$ under $B$, that
\begin{equation}\label{equation3}
 e_1.e_3=-\frac{2}{3}e_1.e_2=e_2.e_4
\end{equation}
We now obtain from (\ref{equation2}) and (\ref{equation3}), with respect to the standard basis $\{e_1,e_2,e_3,e_4\}$, the matrix form of {\scriptsize $$\Omega=\begin{pmatrix}
\begin {array}{rrrr} 0&1&-2/3&-1\\ \noalign{\medskip}-1&0&1&-2
/3\\ \noalign{\medskip}2/3&-1&0&1\\ \noalign{\medskip}1&2/3&-1&0
\end {array}
\end{pmatrix}.$$}

\subsection*{Proof of the arithmeticity of $\Gamma(f,g)$} 
By computation, we obtain that $\epsilon_1=3e_1+e_3$, $\epsilon_2=3e_1+2e_2+3e_3$, $\epsilon_2^*=3e_1+4e_2+3e_3+2e_4$ and $\epsilon_1^*=e_1$ form a basis of $\Q^4$ over $\Q$, with respect to which 
{\scriptsize $$\Omega=\begin{pmatrix}\begin {array}{rrrr} 0&0&0&2/3\\ \noalign{\medskip}0&0&-8/3&0
\\ \noalign{\medskip}0&8/3&0&0\\ \noalign{\medskip}-2/3&0&0&0
\end {array}
\end{pmatrix}.$$}
Let $$P=C=A^{-1}B,\quad  Q=A^{-7}B^3CB^{-3}A^7,\quad R=B^{-3}A^7CA^{-7}B^3.$$ It can be
checked     easily    that,     with    respect     to     the    basis
$\{\epsilon_1,\epsilon_2,\epsilon_2^*,\epsilon_1^*\}$,  the $P, Q, R$ have, respectively, the matrix form
{\scriptsize \[\begin{pmatrix} \begin {array}{rrrr} 1&0&0&0\\ \noalign{\medskip}0&1&2&0
\\ \noalign{\medskip}0&0&1&0\\ \noalign{\medskip}0&0&0&1\end {array}
 \end{pmatrix}, \qquad  \begin{pmatrix} \begin {array}{rrrr} 1&0&0&0\\ \noalign{\medskip}0&1&0&0
\\ \noalign{\medskip}0&-2&1&0\\ \noalign{\medskip}0&0&0&1\end {array}
   \end{pmatrix}, \qquad \begin{pmatrix}  \begin {array}{rrrr} 1&-24&288&-72\\ \noalign{\medskip}0&-23&
288&-72\\ \noalign{\medskip}0&-2&25&-6\\ \noalign{\medskip}0&0&0&1
\end {array}
 \end{pmatrix}.\]}  A  computation  shows  that  if
\[E=Q^{-1}P^{-6}RP^6, \quad F=PEP^{-1},\] \[x=[E,F]=EFE^{-1}F^{-1},\quad y=E^8x,\quad u=F^8y^{-1}, \quad z=u^{-576}x^{17856},\] then

{\scriptsize\[E=\begin{pmatrix}
\begin {array}{rrrr} 1&-24&0&-72\\ \noalign{\medskip}0&1&0&0
\\ \noalign{\medskip}0&0&1&-6\\ \noalign{\medskip}0&0&0&1\end {array}
\end{pmatrix},\quad F=\begin{pmatrix}
\begin {array}{rrrr} 1&-24&48&-72\\ \noalign{\medskip}0&1&0&-
12\\ \noalign{\medskip}0&0&1&-6\\ \noalign{\medskip}0&0&0&1
\end {array}
\end{pmatrix},\]} {\scriptsize\[x=\begin{pmatrix} \begin {array}{rrrr} 1&0&0&576\\ \noalign{\medskip}0&1&0&0
\\ \noalign{\medskip}0&0&1&0\\ \noalign{\medskip}0&0&0&1\end {array}
\end{pmatrix}, \quad y=\begin{pmatrix} \begin {array}{rrrr} 1&-192&0&0\\ \noalign{\medskip}0&1&0&0
\\ \noalign{\medskip}0&0&1&-48\\ \noalign{\medskip}0&0&0&1\end {array}
\end{pmatrix},\]} {\scriptsize\[u=\begin{pmatrix}
             \begin {array}{rrrr} 1&0&384&17856\\ \noalign{\medskip}0&1&0&-
96\\ \noalign{\medskip}0&0&1&0\\ \noalign{\medskip}0&0&0&1\end {array}
            \end{pmatrix},\quad z=\begin{pmatrix} \begin {array}{rrrr} 1&0&-221184&0\\ \noalign{\medskip}0&1&0&
55296\\ \noalign{\medskip}0&0&1&0\\ \noalign{\medskip}0&0&0&1
\end {array}
\end{pmatrix}.\]} Observe that $P, x, y$ and $z$ are (non-trivial) unipotent elements in $\Gamma(f,g)$, which correspond to the positive roots of $\Sp_4(\Omega)$; and among them $x$ and $z$  correspond, respectively, to the highest and second highest roots. Since $\Gamma(f,g)$ is Zariski dense in
$\Sp_4(\Omega)$ by \cite{BH}, it follows  from \cite{T} (cf. \cite[Theorem 3.5]{Ve}, Remark \ref{Venkataramana})
that $\Gamma(f,g)$ is an arithmetic subgroup of $\Sp_4(\Omega)(\Z)$. \qed

\subsection{Example \ref{arithmetic13} of Table \ref{table:newarithmetic}}\label{proofarithmetic13}
This is Example 27 of \cite[Table 2]{SV}. In this case  $$\alpha=\left(\frac{1}{3},\frac{1}{3},\frac{2}{3},\frac{2}{3}\right), \quad \beta=\left(0,0,\frac{1}{4},\frac{3}{4}\right);$$  $$f(X)=X^4+2X^3+3X^2+2X+1, \quad g(X)=X^4-2X^3+2X^2-2X+1;$$ and $f(X)-g(X)=4X^3+X^2+4X$.

Let $A$ and $B$  be the companion  matrices of $f(X)$ and  $g(X)$ respectively, and let $C=A^{-1}B$. Then
{\scriptsize \[A=\begin{pmatrix}\begin {array}{rrrr} 0&0&0&-1\\ \noalign{\medskip}1&0&0&-2
\\ \noalign{\medskip}0&1&0&-3\\ \noalign{\medskip}0&0&1&-2\end {array}
  \end{pmatrix},\quad  B=\begin{pmatrix}\begin {array}{rrrr} 0&0&0&-1\\ \noalign{\medskip}1&0&0&2
\\ \noalign{\medskip}0&1&0&-2\\ \noalign{\medskip}0&0&1&2\end {array}
 \end{pmatrix},\quad C=\begin{pmatrix}\begin {array}{rrrr} 1&0&0&4\\ \noalign{\medskip}0&1&0&1
\\ \noalign{\medskip}0&0&1&4\\ \noalign{\medskip}0&0&0&1\end {array}
 \end{pmatrix}.\]}   Let $\Gamma(f,g)$ be the  subgroup  of
$\SL_4(\Z)$ generated by $A$ and $B$. 

\subsection*{The invariant symplectic form} Using the same method as in Subsection \ref{proofarithmetic18}, with respect to the standard basis $\{e_1,e_2,e_3,e_4\}$, we obtain the matrix form of
{\scriptsize\[\Omega=\begin{pmatrix}
 \begin {array}{rrrr} 0&-4&1&6\\ \noalign{\medskip}4&0&-4&1
\\ \noalign{\medskip}-1&4&0&-4\\ \noalign{\medskip}-6&-1&4&0
\end {array}
\end{pmatrix}.\]}
\subsection*{Proof of the arithmeticity of $\Gamma(f,g)$}
By computation, we obtain that $\epsilon_1=e_1+e_3$, $\epsilon_2=4e_1+e_2+4e_3$, $\epsilon_2^*=-10e_1-8e_2-8e_3-e_4$ and $\epsilon_1^*=e_1+2e_2$ form a basis of $\Q^4$ over $\Q$, with respect to which, the matrix form of  
{\scriptsize $$\Omega=\begin{pmatrix}\begin {array}{rrrr} 0&0&0&-1\\ \noalign{\medskip}0&0&-9&0
\\ \noalign{\medskip}0&9&0&0\\ \noalign{\medskip}1&0&0&0\end {array}
 \end{pmatrix}.$$}
Let $$P=C=A^{-1}B,\quad  Q=A^{-4}CA^4,\quad R=A^4CA^{-4}.$$ It can be
checked     easily    that     with    respect     to     the    basis
$\{\epsilon_1,\epsilon_2,\epsilon_2^*,\epsilon_1^*\}$,  the $P, Q, R$ have, respectively, the matrix form
{\tiny \[\begin{pmatrix}\begin {array}{rrrr} 1&0&0&0\\ \noalign{\medskip}0&1&-1&0
\\ \noalign{\medskip}0&0&1&0\\ \noalign{\medskip}0&0&0&1\end {array}
 \end{pmatrix},\quad  \begin{pmatrix} \begin {array}{rrrr} 1&0&0&0\\ \noalign{\medskip}0&1&0&0
\\ \noalign{\medskip}0&1&1&0\\ \noalign{\medskip}0&0&0&1\end {array}
   \end{pmatrix},\quad \begin{pmatrix}\begin {array}{rrrr} 1&-36&468&-144\\ \noalign{\medskip}0&14&-
169&52\\ \noalign{\medskip}0&1&-12&4\\ \noalign{\medskip}0&0&0&1
\end {array}
 \end{pmatrix}.\]}  A  computation  shows  that  if
\[S=P^{13}RP^{-13}Q^{-1},\ T=PSP^{-1},\ x=[T,S],\ y=S^{-2}x,\ z=T^2yx,\] then
{\scriptsize\[S=\begin{pmatrix}\begin {array}{rrrr} 1&-36&0&-144\\ \noalign{\medskip}0&1&0&0
\\ \noalign{\medskip}0&0&1&4\\ \noalign{\medskip}0&0&0&1\end {array}
\end{pmatrix},\ T=\begin{pmatrix}
 \begin {array}{rrrr} 1&-36&-36&-144\\ \noalign{\medskip}0&1&0&
-4\\ \noalign{\medskip}0&0&1&4\\ \noalign{\medskip}0&0&0&1\end {array}
\end{pmatrix},\]} {\scriptsize\[x=\begin{pmatrix}
\begin {array}{rrrr} 1&0&0&-288\\ \noalign{\medskip}0&1&0&0
\\ \noalign{\medskip}0&0&1&0\\ \noalign{\medskip}0&0&0&1\end {array}
\end{pmatrix},\ y=\begin{pmatrix}\begin {array}{rrrr} 1&72&0&0\\ \noalign{\medskip}0&1&0&0
\\ \noalign{\medskip}0&0&1&-8\\ \noalign{\medskip}0&0&0&1\end {array}
\end{pmatrix},\ z=\begin{pmatrix}\begin {array}{rrrr} 1&0&-72&0\\ \noalign{\medskip}0&1&0&-8
\\ \noalign{\medskip}0&0&1&0\\ \noalign{\medskip}0&0&0&1\end {array}
\end{pmatrix}.\]} Observe that $P, x, y$ and $z$ are (non-trivial) unipotent elements in $\Gamma(f,g)$, which correspond to the positive roots of $\Sp_4(\Omega)$; and among them $x$ and $z$  correspond, respectively, to the highest and second highest roots. Since $\Gamma(f,g)$ is Zariski dense in
$\Sp_4(\Omega)$ by \cite{BH}, it follows  from \cite{T} (cf. \cite[Theorem 3.5]{Ve}, Remark \ref{Venkataramana})
that $\Gamma(f,g)$ is an arithmetic subgroup of $\Sp_4(\Omega)(\Z)$. \qed

\subsection{Example \ref{arithmetic14} of Table \ref{table:newarithmetic}}\label{proofarithmetic14}
This is Example 28 of \cite[Table 2]{SV}. In this case  $$\alpha=\left(\frac{1}{3},\frac{1}{3},\frac{2}{3},\frac{2}{3}\right), \quad \beta=\left(0,0,\frac{1}{6},\frac{5}{6}\right);$$  $$f(X)=X^4+2X^3+3X^2+2X+1, \quad g(X)=X^4-3X^3+4X^2-3X+1;$$ and $f(X)-g(X)=5X^3-X^2+5X$.

Let $A$ and $B$  be the companion  matrices of $f(X)$ and  $g(X)$ respectively, and let $C=A^{-1}B$. Then
{\scriptsize \[A=\begin{pmatrix}\begin {array}{rrrr} 0&0&0&-1\\ \noalign{\medskip}1&0&0&-2
\\ \noalign{\medskip}0&1&0&-3\\ \noalign{\medskip}0&0&1&-2\end {array}
  \end{pmatrix},\quad  B=\begin{pmatrix}\begin {array}{rrrr} 0&0&0&-1\\ \noalign{\medskip}1&0&0&3
\\ \noalign{\medskip}0&1&0&-4\\ \noalign{\medskip}0&0&1&3\end {array}
 \end{pmatrix},\quad C=\begin{pmatrix}\begin {array}{rrrr} 1&0&0&5\\ \noalign{\medskip}0&1&0&-1
\\ \noalign{\medskip}0&0&1&5\\ \noalign{\medskip}0&0&0&1\end {array}
 \end{pmatrix}.\]}   Let $\Gamma(f,g)$ be the  subgroup  of
$\SL_4(\Z)$ generated by $A$ and $B$. 

\subsection*{The invariant symplectic form} Using the same method as in Subsection \ref{proofarithmetic18},  with respect to the standard basis $\{e_1,e_2,e_3,e_4\}$, we obtain the matrix form of
{\scriptsize\[\Omega=\begin{pmatrix}
 \begin {array}{rrrr} 0&5&1&-12\\ \noalign{\medskip}-5&0&5&1
\\ \noalign{\medskip}-1&-5&0&5\\ \noalign{\medskip}12&-1&-5&0
\end {array}
\end{pmatrix}.\]}
\subsection*{Proof of the arithmeticity of $\Gamma(f,g)$}
By computation, we obtain that $\epsilon_1=e_1+e_3$, $\epsilon_2=5e_1-e_2+5e_3$, $\epsilon_2^*=-17e_1-10e_2-10e_3+e_4$ and $\epsilon_1^*=e_1+2e_2$ form a basis of $\Q^4$ over $\Q$, with respect to which, the matrix form of  
{\scriptsize $$\Omega=\begin{pmatrix}\begin {array}{rrrr} 0&0&0&-1\\ \noalign{\medskip}0&0&-36&0
\\ \noalign{\medskip}0&36&0&0\\ \noalign{\medskip}1&0&0&0\end {array}
 \end{pmatrix}.$$}
Let $$P=C=A^{-1}B,\quad  Q=A^{-4}CA^4,\quad R=A^4CA^{-4}.$$   It can be
checked     easily    that     with    respect     to     the    basis
$\{\epsilon_1,\epsilon_2,\epsilon_2^*,\epsilon_1^*\}$,  the $P, Q, R$ have, respectively, the matrix form
{\tiny \[\begin{pmatrix}\begin {array}{rrrr} 1&0&0&0\\ \noalign{\medskip}0&1&1&0
\\ \noalign{\medskip}0&0&1&0\\ \noalign{\medskip}0&0&0&1\end {array}
 \end{pmatrix},\quad  \begin{pmatrix} \begin {array}{rrrr} 1&0&0&0\\ \noalign{\medskip}0&1&0&0
\\ \noalign{\medskip}0&-1&1&0\\ \noalign{\medskip}0&0&0&1\end {array}
\end{pmatrix},\quad \begin{pmatrix}\begin {array}{rrrr} 1&-144&-3312&576\\ \noalign{\medskip}0&24
&529&-92\\ \noalign{\medskip}0&-1&-22&4\\ \noalign{\medskip}0&0&0&1
\end {array}
 \end{pmatrix}.\]}  A  computation  shows  that  if
\[S=P^{23}RP^{-23}Q^{-1},\ T=PSP^{-1}, \ x=[T,S],\ y=S^{-2}x, \ z=T^2yx,\] then
{\scriptsize\[S=\begin{pmatrix}\begin {array}{rrrr} 1&-144&0&576\\ \noalign{\medskip}0&1&0&0
\\ \noalign{\medskip}0&0&1&4\\ \noalign{\medskip}0&0&0&1\end {array}
\end{pmatrix},\quad T=\begin{pmatrix}
 \begin {array}{rrrr} 1&-144&144&576\\ \noalign{\medskip}0&1&0&
4\\ \noalign{\medskip}0&0&1&4\\ \noalign{\medskip}0&0&0&1\end {array}
\end{pmatrix},\]} {\scriptsize\[x=\begin{pmatrix}
\begin {array}{rrrr} 1&0&0&1152\\ \noalign{\medskip}0&1&0&0
\\ \noalign{\medskip}0&0&1&0\\ \noalign{\medskip}0&0&0&1\end {array}
\end{pmatrix},\ y=\begin{pmatrix}\begin {array}{rrrr} 1&288&0&0\\ \noalign{\medskip}0&1&0&0
\\ \noalign{\medskip}0&0&1&-8\\ \noalign{\medskip}0&0&0&1\end {array}
\end{pmatrix},\ z=\begin{pmatrix}\begin {array}{rrrr} 1&0&288&0\\ \noalign{\medskip}0&1&0&8
\\ \noalign{\medskip}0&0&1&0\\ \noalign{\medskip}0&0&0&1\end {array}
\end{pmatrix}.\]} Observe that $P, x, y$ and $z$ are (non-trivial) unipotent elements in $\Gamma(f,g)$, which correspond to the positive roots of $\Sp_4(\Omega)$; and among them $x$ and $z$  correspond, respectively, to the highest and second highest roots. Since $\Gamma(f,g)$ is Zariski dense in
$\Sp_4(\Omega)$ by \cite{BH}, it follows  from \cite{T} (cf. \cite[Theorem 3.5]{Ve}, Remark \ref{Venkataramana})
that $\Gamma(f,g)$ is an arithmetic subgroup of $\Sp_4(\Omega)(\Z)$. \qed

\subsection{Example \ref{arithmetic15} of Table \ref{table:newarithmetic}}\label{proofarithmetic15}
This is Example 29 of \cite[Table 2]{SV}. In this case  $$\alpha=\left(\frac{1}{3},\frac{1}{3},\frac{2}{3},\frac{2}{3}\right), \quad \beta=\left(\frac{1}{6},\frac{1}{6},\frac{5}{6},\frac{5}{6}\right);$$  $$f(X)=X^4+2X^3+3X^2+2X+1, \quad g(X)=X^4-2X^3+3X^2-2X+1;$$ and $f(X)-g(X)=4X^3+4X$.

Let $A$ and $B$  be the companion  matrices of $f(X)$ and  $g(X)$ respectively, and let $C=A^{-1}B$. Then
{\scriptsize \[A=\begin{pmatrix}\begin {array}{rrrr} 0&0&0&-1\\ \noalign{\medskip}1&0&0&-2
\\ \noalign{\medskip}0&1&0&-3\\ \noalign{\medskip}0&0&1&-2\end {array}
  \end{pmatrix},\quad  B=\begin{pmatrix}\begin {array}{rrrr} 0&0&0&-1\\ \noalign{\medskip}1&0&0&2
\\ \noalign{\medskip}0&1&0&-3\\ \noalign{\medskip}0&0&1&2\end {array}
 \end{pmatrix},\quad C=\begin{pmatrix}\begin {array}{rrrr} 1&0&0&4\\ \noalign{\medskip}0&1&0&0
\\ \noalign{\medskip}0&0&1&4\\ \noalign{\medskip}0&0&0&1\end {array}
 \end{pmatrix}.\]}   Let $\Gamma(f,g)$ be the  subgroup  of
$\SL_4(\Z)$ generated by $A$ and $B$. 

\subsection*{The invariant symplectic form} Using the same method as in Subsection \ref{proofarithmetic18},  with respect to the standard basis $\{e_1,e_2,e_3,e_4\}$, we obtain the matrix form of
{\scriptsize\[\Omega=\begin{pmatrix}
         \begin {array}{rrrr} 0&-1/2&0&1\\ \noalign{\medskip}1/2&0&-1/2
&0\\ \noalign{\medskip}0&1/2&0&-1/2\\ \noalign{\medskip}-1&0&1/2&0
\end {array}
\end{pmatrix}.\]}
\subsection*{Proof of the arithmeticity of $\Gamma(f,g)$}
By computation, we obtain that $\epsilon_1=2e_1+e_2+2e_3$, $\epsilon_2=4e_1+4e_3$, $\epsilon_2^*=-16e_1-8e_2-8e_3+4e_4$ and $\epsilon_1^*=e_1+2e_2$ form a basis of $\Q^4$ over $\Q$, with respect to which, the matrix form of 
{\scriptsize $$\Omega=\begin{pmatrix}\begin {array}{rrrr} 0&0&0&1/2\\ \noalign{\medskip}0&0&8&0
\\ \noalign{\medskip}0&-8&0&0\\ \noalign{\medskip}-1/2&0&0&0
\end {array}
\end{pmatrix}.$$}
Let $$P=C=A^{-1}B,\quad  Q=A^{-7}CA^7,\quad R=B^{-1}CB.$$   It can be
checked     easily    that     with    respect     to     the    basis
$\{\epsilon_1,\epsilon_2,\epsilon_2^*,\epsilon_1^*\}$,  the $P, Q, R$ have, respectively, the matrix form
{\scriptsize \[\begin{pmatrix} \begin {array}{rrrr} 1&0&0&0\\ \noalign{\medskip}0&1&4&0
\\ \noalign{\medskip}0&0&1&0\\ \noalign{\medskip}0&0&0&1\end {array}
\end{pmatrix},\quad \begin{pmatrix}\begin {array}{rrrr} 1&0&0&0\\ \noalign{\medskip}0&1&0&0
\\ \noalign{\medskip}0&-4&1&0\\ \noalign{\medskip}0&0&0&1\end {array}
   \end{pmatrix},\quad \begin{pmatrix}\begin {array}{rrrr} 1&0&0&0\\ \noalign{\medskip}0&1&0&0
\\ \noalign{\medskip}-2&-4&1&0\\ \noalign{\medskip}-16&-32&0&1
\end {array}
\end{pmatrix}.\]}  A  computation  shows  that  if
\[S=RQ^{-1},\ T=PSP^{-1},\ x=[T,S],\ y=S^{32}x^{-1},\ z=x^{511}T^{32}y^{-1},\] then
{\scriptsize\[S= \begin{pmatrix}\begin {array}{rrrr} 1&0&0&0\\ \noalign{\medskip}0&1&0&0
\\ \noalign{\medskip}-2&0&1&0\\ \noalign{\medskip}-16&-32&0&1
\end {array}
\end{pmatrix},\quad T=\begin{pmatrix}\begin {array}{rrrr} 1&0&0&0\\ \noalign{\medskip}-8&1&0&0
\\ \noalign{\medskip}-2&0&1&0\\ \noalign{\medskip}-16&-32&128&1
\end {array}
\end{pmatrix},\]} {\scriptsize\[x=\begin{pmatrix}\begin {array}{rrrr} 1&0&0&0\\ \noalign{\medskip}0&1&0&0
\\ \noalign{\medskip}0&0&1&0\\ \noalign{\medskip}-512&0&0&1
\end {array}
\end{pmatrix}, y=\begin{pmatrix} \begin {array}{rrrr} 1&0&0&0\\ \noalign{\medskip}0&1&0&0
\\ \noalign{\medskip}-64&0&1&0\\ \noalign{\medskip}0&-1024&0&1
\end {array}
\end{pmatrix},  z=\begin{pmatrix} \begin {array}{rrrr} 1&0&0&0\\ \noalign{\medskip}-256&1&0&0
\\ \noalign{\medskip}0&0&1&0\\ \noalign{\medskip}0&0&4096&1
\end {array}
\end{pmatrix}.\]} 

In this case, note that if we change the basis $\{\epsilon_1,\epsilon_2,\epsilon_2^*,\epsilon_1^*\}$ to a new basis $\{\epsilon_1^*,\epsilon_2^*,\epsilon_2,\epsilon_1\}$ (that is, we take the transpose of all the above matrices), then the matrix form of $\Omega$, with respect to the new basis, is still anti-diagonal; and $Q, x, y, z$, with respect to the new basis, are (non-trivial) unipotent elements in $\Gamma(f,g)$, which correspond to the positive roots of $\Sp_4(\Omega)$; and among them $x$ and $y$  correspond, respectively, to the highest and second highest roots. Since $\Gamma(f,g)$ is Zariski dense in
$\Sp_4(\Omega)$ by \cite{BH}, it follows  from \cite{T} (cf. \cite[Theorem 3.5]{Ve}, Remark \ref{Venkataramana})
that $\Gamma(f,g)$ is an arithmetic subgroup of $\Sp_4(\Omega)(\Z)$. \qed

\subsection{Example \ref{arithmetic16} of Table \ref{table:newarithmetic}}\label{proofarithmetic16}
This is Example 30 of \cite[Table 2]{SV}. In this case  $$\alpha=\left(\frac{1}{3},\frac{1}{3},\frac{2}{3},\frac{2}{3}\right), \quad \beta=\left(\frac{1}{4},\frac{3}{4},\frac{1}{6},\frac{5}{6}\right);$$  $$f(X)=X^4+2X^3+3X^2+2X+1, \quad g(X)=X^4-X^3+2X^2-X+1;$$ and $f(X)-g(X)=3X^3+X^2+3X$.

Let $A$ and $B$  be the companion  matrices of $f(X)$ and  $g(X)$ respectively, and let $C=A^{-1}B$. Then
{\scriptsize \[A=\begin{pmatrix}\begin {array}{rrrr} 0&0&0&-1\\ \noalign{\medskip}1&0&0&-2
\\ \noalign{\medskip}0&1&0&-3\\ \noalign{\medskip}0&0&1&-2\end {array}
  \end{pmatrix},\  B=\begin{pmatrix}\begin {array}{rrrr} 0&0&0&-1\\ \noalign{\medskip}1&0&0&1
\\ \noalign{\medskip}0&1&0&-2\\ \noalign{\medskip}0&0&1&1\end {array}
 \end{pmatrix},\ C=\begin{pmatrix}\begin {array}{rrrr} 1&0&0&3\\ \noalign{\medskip}0&1&0&1
\\ \noalign{\medskip}0&0&1&3\\ \noalign{\medskip}0&0&0&1\end {array}
 \end{pmatrix}.\]}   Let $\Gamma(f,g)$ be the  subgroup  of
$\SL_4(\Z)$ generated by $A$ and $B$. 

\subsection*{The invariant symplectic form} Using the same method as in Subsection \ref{proofarithmetic18}, with respect to the standard basis $\{e_1,e_2,e_3,e_4\}$, we obtain the matrix form of
{\scriptsize\[\Omega=\begin{pmatrix}
          \begin {array}{rrrr} 0&-3/4&1/4&1\\ \noalign{\medskip}3/4&0&-3
/4&1/4\\ \noalign{\medskip}-1/4&3/4&0&-3/4\\ \noalign{\medskip}-1&-1/4
&3/4&0\end {array}
\end{pmatrix}.\]}
\subsection*{Proof of the arithmeticity of $\Gamma(f,g)$}
By computation, we obtain that $\epsilon_1=e_1+e_3$, $\epsilon_2=3e_1+e_2+3e_3$, $\epsilon_2^*=-4e_1-e_2-3e_3-e_4$ and $\epsilon_1^*=e_1-e_2$ form a basis of $\Q^4$ over $\Q$, with respect to which, the matrix form of 
{\scriptsize $$\Omega=\begin{pmatrix}\begin {array}{rrrr} 0&0&0&-1/4\\ \noalign{\medskip}0&0&-1&0
\\ \noalign{\medskip}0&1&0&0\\ \noalign{\medskip}1/4&0&0&0\end {array}
\end{pmatrix}.$$}
Let $$P=C=A^{-1}B,\quad  Q=B^{-3}CB^3,\quad R=B^3CB^{-3}.$$   It can be
checked     easily    that     with    respect     to     the    basis
$\{\epsilon_1,\epsilon_2,\epsilon_2^*,\epsilon_1^*\}$,  the $P, Q, R$ have, respectively, the matrix form
{\scriptsize \[\begin{pmatrix} \begin {array}{rrrr} 1&0&0&0\\ \noalign{\medskip}0&1&-1&0
\\ \noalign{\medskip}0&0&1&0\\ \noalign{\medskip}0&0&0&1\end {array}
 \end{pmatrix}, \qquad  \begin{pmatrix} \begin {array}{rrrr} 1&0&0&0\\ \noalign{\medskip}0&1&0&0
\\ \noalign{\medskip}0&1&1&0\\ \noalign{\medskip}0&0&0&1\end {array}
   \end{pmatrix}, \qquad \begin{pmatrix}\begin {array}{rrrr} 1&8&0&-16\\ \noalign{\medskip}0&1&0&0
\\ \noalign{\medskip}0&1&1&-2\\ \noalign{\medskip}0&0&0&1\end {array}
 \end{pmatrix}.\]}  A  computation  shows  that  if
\[E=PQ^{-1}RP^{-1}, \ F=RQ^{-1}, \ x=[E,F],\ y=F^{-2}x, \ z=E^2yx,\] then
{\scriptsize\[E=\begin{pmatrix}
\begin {array}{rrrr} 1&8&8&-16\\ \noalign{\medskip}0&1&0&2
\\ \noalign{\medskip}0&0&1&-2\\ \noalign{\medskip}0&0&0&1\end {array}
\end{pmatrix},\qquad F=\begin{pmatrix}
\begin {array}{rrrr} 1&8&0&-16\\ \noalign{\medskip}0&1&0&0
\\ \noalign{\medskip}0&0&1&-2\\ \noalign{\medskip}0&0&0&1\end {array}
\end{pmatrix},\]} {\scriptsize\[x=\begin{pmatrix} \begin {array}{rrrr} 1&0&0&-32\\ \noalign{\medskip}0&1&0&0
\\ \noalign{\medskip}0&0&1&0\\ \noalign{\medskip}0&0&0&1\end {array}
\end{pmatrix}, \quad y=\begin{pmatrix} \begin {array}{rrrr} 1&-16&0&0\\ \noalign{\medskip}0&1&0&0
\\ \noalign{\medskip}0&0&1&4\\ \noalign{\medskip}0&0&0&1\end {array}
\end{pmatrix},\quad z=\begin{pmatrix} \begin {array}{rrrr} 1&0&16&0\\ \noalign{\medskip}0&1&0&4
\\ \noalign{\medskip}0&0&1&0\\ \noalign{\medskip}0&0&0&1\end {array}
\end{pmatrix}.\]} Observe that $P, x, y$ and $z$ are (non-trivial) unipotent elements in $\Gamma(f,g)$, which correspond to the positive roots of $\Sp_4(\Omega)$; and among them $x$ and $z$  correspond, respectively, to the highest and second highest roots. Since $\Gamma(f,g)$ is Zariski dense in
$\Sp_4(\Omega)$ by \cite{BH}, it follows  from \cite{T} (cf. \cite[Theorem 3.5]{Ve}, Remark \ref{Venkataramana})
that $\Gamma(f,g)$ is an arithmetic subgroup of $\Sp_4(\Omega)(\Z)$. \qed

\subsection{Example \ref{arithmetic17} of Table \ref{table:newarithmetic}}\label{proofarithmetic17}
This is Example 31 of \cite[Table 2]{SV}. In this case  $$\alpha=\left(\frac{1}{3},\frac{1}{3},\frac{2}{3},\frac{2}{3}\right), \quad \beta=\left(\frac{1}{10},\frac{3}{10},\frac{7}{10},\frac{9}{10}\right);$$  $$f(X)=X^4+2X^3+3X^2+2X+1, \quad g(X)=X^4-X^3+X^2-X+1;$$ and $f(X)-g(X)=3X^3+2X^2+3X$.\\

Let $A$ and $B$  be the companion  matrices of $f(X)$ and  $g(X)$ respectively, and let $C=A^{-1}B$. Then
{\scriptsize \[A=\begin{pmatrix}\begin {array}{rrrr} 0&0&0&-1\\ \noalign{\medskip}1&0&0&-2
\\ \noalign{\medskip}0&1&0&-3\\ \noalign{\medskip}0&0&1&-2\end {array}
  \end{pmatrix},\  B=\begin{pmatrix}\begin {array}{rrrr} 0&0&0&-1\\ \noalign{\medskip}1&0&0&1
\\ \noalign{\medskip}0&1&0&-1\\ \noalign{\medskip}0&0&1&1\end {array}
 \end{pmatrix},\ C=\begin{pmatrix}\begin {array}{rrrr} 1&0&0&3\\ \noalign{\medskip}0&1&0&2
\\ \noalign{\medskip}0&0&1&3\\ \noalign{\medskip}0&0&0&1\end {array}
 \end{pmatrix}.\]}   Let $\Gamma(f,g)$ be the  subgroup  of
$\SL_4(\Z)$ generated by $A$ and $B$. 

\subsection*{The invariant symplectic form} Using the same method as in Subsection \ref{proofarithmetic18}, we obtain the matrix form of
{\scriptsize\[\Omega=\begin{pmatrix}
          \begin {array}{rrrr} 0&-3/2&1&1\\ \noalign{\medskip}3/2&0&-3/2
&1\\ \noalign{\medskip}-1&3/2&0&-3/2\\ \noalign{\medskip}-1&-1&3/2&0
\end {array}
\end{pmatrix}.\]}
\subsection*{Proof of the arithmeticity of $\Gamma(f,g)$}
By computation, we obtain that $\epsilon_1=e_1+e_3$, $\epsilon_2=3e_1+2e_2+3e_3$, $\epsilon_2^*=-5e_1-6e_2-6e_3-2e_4$ and $\epsilon_1^*=e_1+2e_2$ form a basis of $\Q^4$ over $\Q$, with respect to which, the matrix form of 
{\scriptsize $$\Omega=\begin{pmatrix}\begin {array}{rrrr} 0&0&0&-1\\ \noalign{\medskip}0&0&-1&0
\\ \noalign{\medskip}0&1&0&0\\ \noalign{\medskip}1&0&0&0\end {array}
\end{pmatrix}.$$}
Let $$P=C=A^{-1}B,\quad  Q=A^{-4}CA^4,\quad R=A^4CA^{-4},\quad S=P^2RP^{-1}Q.$$   It can be
checked     easily    that     with    respect     to     the    basis
$\{\epsilon_1,\epsilon_2,\epsilon_2^*,\epsilon_1^*\}$,  the $P, Q, R, S$ have, respectively, the matrix form
{\tiny \[\begin{pmatrix} \begin {array}{rrrr} 1&0&0&0\\ \noalign{\medskip}0&1&-2&0
\\ \noalign{\medskip}0&0&1&0\\ \noalign{\medskip}0&0&0&1\end {array}
 \end{pmatrix},  \begin{pmatrix} \begin {array}{rrrr} 1&0&0&0\\ \noalign{\medskip}0&1&0&0
\\ \noalign{\medskip}0&2&1&0\\ \noalign{\medskip}0&0&0&1\end {array}
   \end{pmatrix}, \begin{pmatrix}\begin {array}{rrrr} 1&-4&12&-8\\ \noalign{\medskip}0&7&-18&12
\\ \noalign{\medskip}0&2&-5&4\\ \noalign{\medskip}0&0&0&1\end {array}
 \end{pmatrix}, \begin{pmatrix}\begin {array}{rrrr} 1&4&4&-8\\ \noalign{\medskip}0&-1&0&-4
\\ \noalign{\medskip}0&0&-1&4\\ \noalign{\medskip}0&0&0&1\end {array}
\end{pmatrix}.\]}  A  computation  shows  that  if
\[E=[R,S], \quad F=QEQ^{-1}, \quad x=[E,F],\quad y=(E^{-1}F)^2x,\] \[u=F^{-288}y^{168}, \quad z=u^{20736}x^{1003401216},\] then
{\scriptsize\[E=\begin{pmatrix}
\begin {array}{rrrr} 1&24&-72&384\\ \noalign{\medskip}0&1&0&-
72\\ \noalign{\medskip}0&0&1&-24\\ \noalign{\medskip}0&0&0&1\end{array}
\end{pmatrix},\qquad F=\begin{pmatrix}
\begin {array}{rrrr} 1&168&-72&384\\ \noalign{\medskip}0&1&0&-
72\\ \noalign{\medskip}0&0&1&-168\\ \noalign{\medskip}0&0&0&1\end{array}
\end{pmatrix},\]} {\scriptsize\[ x=\begin{pmatrix}\begin {array}{rrrr} 1&0&0&20736\\ \noalign{\medskip}0&1&0&0
\\ \noalign{\medskip}0&0&1&0\\ \noalign{\medskip}0&0&0&1\end {array}
\end{pmatrix},\qquad y=\begin{pmatrix} \begin {array}{rrrr} 1&288&0&0\\ \noalign{\medskip}0&1&0&0
\\ \noalign{\medskip}0&0&1&-288\\ \noalign{\medskip}0&0&0&1 \end{array}
\end{pmatrix},\]} {\scriptsize\[u=\begin{pmatrix}\begin {array}{rrrr} 1&0&20736&-1003401216
\\ \noalign{\medskip}0&1&0&20736\\ \noalign{\medskip}0&0&1&0
\\ \noalign{\medskip}0&0&0&1\end {array}
\end{pmatrix},\qquad z=\begin{pmatrix}\begin {array}{rrrr} 1&0&429981696&0\\ \noalign{\medskip}0&1&0
&429981696\\ \noalign{\medskip}0&0&1&0\\ \noalign{\medskip}0&0&0&1\end {array}
\end{pmatrix}.\]} Observe that $P, x, y$ and $z$ are (non-trivial) unipotent elements in $\Gamma(f,g)$, which correspond to the positive roots of $\Sp_4(\Omega)$; and among them $x$ and $z$  correspond, respectively, to the highest and second highest roots. Since $\Gamma(f,g)$ is Zariski dense in
$\Sp_4(\Omega)$ by \cite{BH}, it follows  from \cite{T} (cf. \cite[Theorem 3.5]{Ve}, Remark \ref{Venkataramana})
that $\Gamma(f,g)$ is an arithmetic subgroup of $\Sp_4(\Omega)(\Z)$. \qed

\subsection{Example \ref{arithmetic19} of Table \ref{table:newarithmetic}}\label{proofarithmetic20}
This is Example 36 of \cite[Table 2]{SV}. In this case  $$\alpha=\left(\frac{1}{2},\frac{1}{2},\frac{1}{3},\frac{2}{3}\right), \quad \beta=\left(\frac{1}{4},\frac{3}{4},\frac{1}{6},\frac{5}{6}\right);$$  $$f(X)=X^4+3X^3+4X^2+3X+1, \quad g(X)=X^4-X^3+2X^2-X+1;$$ and $f(X)-g(X)=4X^3+2X^2+4X$.

Let $A$ and $B$  be the companion  matrices of $f(X)$ and  $g(X)$ respectively, and let $C=A^{-1}B$. Then
{\scriptsize \[A=\begin{pmatrix}\begin {array}{rrrr} 0&0&0&-1\\ \noalign{\medskip}1&0&0&-3
\\ \noalign{\medskip}0&1&0&-4\\ \noalign{\medskip}0&0&1&-3\end {array}
  \end{pmatrix},\  B=\begin{pmatrix}\begin {array}{rrrr} 0&0&0&-1\\ \noalign{\medskip}1&0&0&1
\\ \noalign{\medskip}0&1&0&-2\\ \noalign{\medskip}0&0&1&1\end {array}
 \end{pmatrix},\ C=\begin{pmatrix}\begin {array}{rrrr} 1&0&0&4\\ \noalign{\medskip}0&1&0&2
\\ \noalign{\medskip}0&0&1&4\\ \noalign{\medskip}0&0&0&1\end {array}
 \end{pmatrix}.\]}   Let $\Gamma(f,g)$ be the  subgroup  of
$\SL_4(\Z)$ generated by $A$ and $B$. 

\subsection*{The invariant symplectic form} Using the same method as in Subsection \ref{proofarithmetic18}, we obtain the matrix form of
{\scriptsize\[\Omega=\begin{pmatrix}
 \begin {array}{rrrr} 0&-2&1&3\\ \noalign{\medskip}2&0&-2&1
\\ \noalign{\medskip}-1&2&0&-2\\ \noalign{\medskip}-3&-1&2&0
\end {array}
\end{pmatrix}.\]}
\subsection*{Proof of the arithmeticity of $\Gamma(f,g)$}
By an easy computation we get $\epsilon_1=e_1+e_3$, $\epsilon_2=4e_1+2e_2+4e_3$, $\epsilon_2^*=-6e_1-2e_2-4e_3-2e_4$ and $\epsilon_1^*=e_1-e_2$ form a basis of $\Q^4$ over $\Q$, with respect to which, the matrix form of 
{\scriptsize $$\Omega=\begin{pmatrix}\begin {array}{rrrr} 0&0&0&-1\\ \noalign{\medskip}0&0&-12&0
\\ \noalign{\medskip}0&12&0&0\\ \noalign{\medskip}1&0&0&0\end {array}
 \end{pmatrix}.$$}
Let $$P=C=A^{-1}B,\quad  Q=B^{-3}CB^3,\quad R=B^3CB^{-3}.$$   It can be
checked     easily    that     with    respect     to     the    basis
$\{\epsilon_1,\epsilon_2,\epsilon_2^*,\epsilon_1^*\}$,  the $P, Q, R$ have, respectively, the matrix form
{\tiny \[\begin{pmatrix}\begin {array}{rrrr} 1&0&0&0\\ \noalign{\medskip}0&1&-2&0
\\ \noalign{\medskip}0&0&1&0\\ \noalign{\medskip}0&0&0&1\end {array}
 \end{pmatrix},\quad  \begin{pmatrix} \begin {array}{rrrr} 1&0&0&0\\ \noalign{\medskip}0&1&0&0
\\ \noalign{\medskip}0&2&1&0\\ \noalign{\medskip}0&0&0&1\end {array}
\end{pmatrix},\quad \begin{pmatrix}\begin {array}{rrrr} 1&24&0&-24\\ \noalign{\medskip}0&1&0&0
\\ \noalign{\medskip}0&2&1&-2\\ \noalign{\medskip}0&0&0&1\end {array}
 \end{pmatrix}.\]}  A  computation  shows  that  if
\[S=RQ^{-1},\ T=PSP^{-1},\ x=[T,S],\ y=S^8x^{-1},\ z=T^{-8}yx^{33},\] then
{\scriptsize\[S=\begin{pmatrix}\begin {array}{rrrr} 1&24&0&-24\\ \noalign{\medskip}0&1&0&0
\\ \noalign{\medskip}0&0&1&-2\\ \noalign{\medskip}0&0&0&1\end {array}
\end{pmatrix},\quad T=\begin{pmatrix}
 \begin {array}{rrrr} 1&24&48&-24\\ \noalign{\medskip}0&1&0&4
\\ \noalign{\medskip}0&0&1&-2\\ \noalign{\medskip}0&0&0&1\end {array}
\end{pmatrix},\]} {\scriptsize\[x=\begin{pmatrix}
\begin {array}{rrrr} 1&0&0&-192\\ \noalign{\medskip}0&1&0&0
\\ \noalign{\medskip}0&0&1&0\\ \noalign{\medskip}0&0&0&1\end {array}
\end{pmatrix}, y=\begin{pmatrix}\begin {array}{rrrr} 1&192&0&0\\ \noalign{\medskip}0&1&0&0
\\ \noalign{\medskip}0&0&1&-16\\ \noalign{\medskip}0&0&0&1\end {array}
\end{pmatrix}, z=\begin{pmatrix}\begin {array}{rrrr} 1&0&-384&0\\ \noalign{\medskip}0&1&0&-32
\\ \noalign{\medskip}0&0&1&0\\ \noalign{\medskip}0&0&0&1\end {array}
\end{pmatrix}.\]} Observe that $P, x, y$ and $z$ are (non-trivial) unipotent elements in $\Gamma(f,g)$, which correspond to the positive roots of $\Sp_4(\Omega)$; and among them $x$ and $z$  correspond, respectively, to the highest and second highest roots. Since $\Gamma(f,g)$ is Zariski dense in
$\Sp_4(\Omega)$ by \cite{BH}, it follows  from \cite{T} (cf. \cite[Theorem 3.5]{Ve}, Remark \ref{Venkataramana})
that $\Gamma(f,g)$ is an arithmetic subgroup of $\Sp_4(\Omega)(\Z)$. \qed

\subsection{Example \ref{arithmetic20} of Table \ref{table:newarithmetic}}\label{proofarithmetic21}
This is Example 39 of \cite[Table 2]{SV}. In this case  $$\alpha=\left(\frac{1}{2},\frac{1}{2},\frac{1}{3},\frac{2}{3}\right), \quad \beta=\left(\frac{1}{12},\frac{5}{12},\frac{7}{12},\frac{11}{12}\right);$$  $$f(X)=X^4+3X^3+4X^2+3X+1, \quad g(X)=X^4-X^2+1;$$ and $f(X)-g(X)=3X^3+5X^2+3X$.

Let $A$ and $B$  be the companion  matrices of $f(X)$ and  $g(X)$ respectively, and let $C=A^{-1}B$. Then
{\scriptsize \[A=\begin{pmatrix}\begin {array}{rrrr} 0&0&0&-1\\ \noalign{\medskip}1&0&0&-3
\\ \noalign{\medskip}0&1&0&-4\\ \noalign{\medskip}0&0&1&-3\end {array}
  \end{pmatrix},\quad  B=\begin{pmatrix}\begin {array}{rrrr} 0&0&0&-1\\ \noalign{\medskip}1&0&0&0
\\ \noalign{\medskip}0&1&0&1\\ \noalign{\medskip}0&0&1&0\end {array}
 \end{pmatrix},\quad C=\begin{pmatrix}\begin {array}{rrrr} 1&0&0&3\\ \noalign{\medskip}0&1&0&5
\\ \noalign{\medskip}0&0&1&3\\ \noalign{\medskip}0&0&0&1\end {array}
 \end{pmatrix}.\]}   Let $\Gamma(f,g)$ be the  subgroup  of
$\SL_4(\Z)$ generated by $A$ and $B$. 

\subsection*{The invariant symplectic form} Using the same method as in Subsection \ref{proofarithmetic18}, we obtain the matrix form of
{\scriptsize\[\Omega=\begin{pmatrix}
 \begin {array}{rrrr} 0&1&-5/3&2\\ \noalign{\medskip}-1&0&1&-5/
3\\ \noalign{\medskip}5/3&-1&0&1\\ \noalign{\medskip}-2&5/3&-1&0
\end {array}
\end{pmatrix}.\]}
\subsection*{Proof of the arithmeticity of $\Gamma(f,g)$}
By an easy computation we obtain that $\epsilon_1=5e_1+9e_2+5e_3$, $\epsilon_2=3e_1+5e_2+3e_3$, $\epsilon_2^*=-3e_1-7e_2-3e_3+2e_4$ and $\epsilon_1^*=5e_1+3e_2$ form a basis of $\Q^4$ over $\Q$, with respect to which, the matrix form of 
{\scriptsize $$\Omega=\begin{pmatrix}\begin {array}{rrrr} 0&0&0&-10/3\\ \noalign{\medskip}0&0&4/3&0
\\ \noalign{\medskip}0&-4/3&0&0\\ \noalign{\medskip}10/3&0&0&0\end{array} \end{pmatrix}.$$}
Let \[P=C=A^{-1}B,\quad  Q=B^{-1}A^{-3}CA^3B,\quad R=A^3BCB^{-1}A^{-3},\] \[S=P^{-1}QP^2R.\]   It can be
checked     easily    that     with    respect     to     the    basis
$\{\epsilon_1,\epsilon_2,\epsilon_2^*,\epsilon_1^*\}$,  the $P, Q, R, S$ have, respectively, the matrix form
{\tiny \[\begin{pmatrix} \begin {array}{rrrr} 1&0&0&0\\ \noalign{\medskip}0&1&2&0
\\ \noalign{\medskip}0&0&1&0\\ \noalign{\medskip}0&0&0&1\end {array}
 \end{pmatrix},  \begin{pmatrix} \begin {array}{rrrr} 1&0&0&0\\ \noalign{\medskip}0&1&0&0
\\ \noalign{\medskip}0&-2&1&0\\ \noalign{\medskip}0&0&0&1\end {array}
   \end{pmatrix}, \begin{pmatrix}\begin {array}{rrrr} 1&-6&-18&-45\\ \noalign{\medskip}0&7&18&
45\\ \noalign{\medskip}0&-2&-5&-15\\ \noalign{\medskip}0&0&0&1
\end {array}
 \end{pmatrix}, \begin{pmatrix}\begin {array}{rrrr} 1&-6&-18&-45\\ \noalign{\medskip}0&-1&0&-
45\\ \noalign{\medskip}0&0&-1&15\\ \noalign{\medskip}0&0&0&1
\end {array}
\end{pmatrix}.\]}  A  computation  shows  that  if
\[E=[R,S], \quad F=QEQ^{-1}, \quad x=[E,F],\quad y=(E^{-1}F)^2x,\] \[u=E^{12}y^{-1}, \quad z=u^{1296}x^{15552},\] then
{\scriptsize\[E=\begin{pmatrix}
\begin {array}{rrrr} 1&-12&-36&0\\ \noalign{\medskip}0&1&0&90
\\ \noalign{\medskip}0&0&1&-30\\ \noalign{\medskip}0&0&0&1\end {array}
\end{pmatrix},\qquad F=\begin{pmatrix}
\begin {array}{rrrr} 1&-84&-36&0\\ \noalign{\medskip}0&1&0&90
\\ \noalign{\medskip}0&0&1&-210\\ \noalign{\medskip}0&0&0&1
\end {array}
\end{pmatrix},\]} {\scriptsize\[ x=\begin{pmatrix}\begin {array}{rrrr} 1&0&0&12960\\ \noalign{\medskip}0&1&0&0
\\ \noalign{\medskip}0&0&1&0\\ \noalign{\medskip}0&0&0&1\end {array}
\end{pmatrix},\qquad y=\begin{pmatrix}\begin {array}{rrrr} 1&-144&0&0\\ \noalign{\medskip}0&1&0&0
\\ \noalign{\medskip}0&0&1&-360\\ \noalign{\medskip}0&0&0&1
\end{array}
\end{pmatrix},\]} {\scriptsize\[u=\begin{pmatrix}\begin {array}{rrrr} 1&0&-432&-155520\\ \noalign{\medskip}0&1&0
&1080\\ \noalign{\medskip}0&0&1&0\\ \noalign{\medskip}0&0&0&1
\end {array} 
\end{pmatrix},\qquad z=\begin{pmatrix}\begin {array}{rrrr} 1&0&-559872&0\\ \noalign{\medskip}0&1&0&
1399680\\ \noalign{\medskip}0&0&1&0\\ \noalign{\medskip}0&0&0&1
\end {array}
\end{pmatrix}.\]} Observe that $P, x, y$ and $z$ are (non-trivial) unipotent elements in $\Gamma(f,g)$, which correspond to the positive roots of $\Sp_4(\Omega)$; and among them $x$ and $z$  correspond, respectively, to the highest and second highest roots. Since $\Gamma(f,g)$ is Zariski dense in
$\Sp_4(\Omega)$ by \cite{BH}, it follows  from \cite{T} (cf. \cite[Theorem 3.5]{Ve}, Remark \ref{Venkataramana})
that $\Gamma(f,g)$ is an arithmetic subgroup of $\Sp_4(\Omega)(\Z)$. \qed

\section*{Acknowledgements}
I am grateful to Duco van Straten, T. N. Venkataramana, and Wadim Zudilin for their constant encouragement and support. The results mentioned in this article were obtained when I was a postdoctoral fellow at the Max Planck Institute for Mathematics in Bonn; I thank the MPI for the postdoctoral fellowship. I also thank Department of Science \& Technology, India for the INSPIRE Faculty Award, and Maple for the computations.

\end{document}